# Optimizing electric vehicles charging through smart energy allocation and cost-saving

Luca Ambrosino[1]
*Dipartimento di Elettronica e Telecomunicazioni (DET)*
*Politecnico di Torino*
Torino, Italy
luca.ambrosino@polito.it

Giuseppe Calafiore
*Dipartimento di Elettronica e Telecomunicazioni (DET) / College of Engineering and Computer Science (CECS)*
*Politecnico di Torino / Vin University*
Torino, Italy / Hanoi, Vietnam
giuseppe.calafiore@polito.it

Khai Manh Nguyen
*College of Engineering and Computer Science (CECS)*
*Vin University*
Hanoi, Vietnam
21khai.nm@vinuni.edu.vn

Riadh Zorgati
*R&D Department*
*Électricité de France (EDF)*
Paris, France
riadh.zorgati@edf.fr

Doanh Nguyen-Ngoc[2]
*College of Engineering and Computer Science (CECS) and College of Enviromental Intelligence (CEI)*
*Vin University*
Hanoi, Vietnam
doanh.nn@vinuni.edu.vn

Laurent El Ghaoui
*College of Engineering and Computer Science (CECS) and College of Enviromental Intelligence (CEI)*
*Vin University*
Hanoi, Vietnam
laurent.eg@vinuni.edu.vn

*Abstract*— As the global focus on combating environmental pollution intensifies, the transition to sustainable energy sources, particularly in the form of electric vehicles (EVs), has become paramount. This paper addresses the pressing need for Smart Charging for EVs by developing a comprehensive mathematical model aimed at optimizing charging station management. The model aims to efficiently allocate the power from charging sockets to EVs, prioritizing cost minimization and avoiding energy waste. Computational simulations demonstrate the efficacy of the mathematical optimization model, which can unleash its full potential when the number of EVs at the charging station is high.

## I. INTRODUCTION

In light of the intensifying environmental crisis, the imperative to transition towards sustainable energy solutions has become more pressing than ever. Central to this transition is the widespread adoption of electric vehicles (EVs) [1], [2], which offer a promising avenue for curbing pollution and mitigating climate change. However, the effective integration of EVs into our daily lives necessitates innovative approaches to address the challenges posed by their charging infrastructure.

Until now, most of the literature concerning electrical vehicles have focused on issues related to electricity and power grid [3]. Data-driven models to reduce EVs impact on the grid itself have been investigated [4], helping to mitigate the risk of overloading the electrical infrastructure finding ways to shift the charging to off-peak hours [5], ensuring a stable and reliable supply of electricity. Furthermore, considering the limited driving range of EVs, they may need to be recharged frequently, especially across long journeys. Consequently, the convenience of charging becomes a primary concern for EV owners and some studies on charging stations location problem have been conducted trying to help this crucial research field [6], [7], even by applying reinforcement learning techniques [8].

On the other hand, our approach deviates from this dominant trend, instead focusing on the optimization of a single charging station when facing the EVs requirements. To predict the charging duration and energy demand of EVs, a hybrid kernel density estimator that uses both Gaussian- and Diffusion-based approach has been developed [9]. Moreover, while there are two types of interaction between vehicles and the charging sockets of charging stations, namely grid-to-vehicle (G2V) and vehicle-to-grid (V2G) [10], our model focuses on grid-to-vehicle interaction. The V2G interaction, on the other hand, allows for the development of more intricate and complex models, for example involving the use of Particle Swarm Optimization (PSO) [11], but has the serious drawback of damaging the integrity and longevity of EV batteries in the long run [12], making this technology less sustainable. This underscores the need for Smart Charging optimization methods to efficiently manage and charge electric vehicles (EVs).

Unlike traditional charging methods, Smart Charging leverages advanced technologies and data-driven algorithms to optimize the utilization of charging stations. By intelligently allocating resources and dynamically adjusting charging schedules, Smart Charging aims to enhance operational efficiency, reduce costs, and alleviate strain on the electrical grid. In essence, Smart Charging transcends the mere act of replenishing the energy reserves of EVs; it represents a holistic approach towards reimagining the entire charging ecosystem. From minimizing peak load demand and optimizing energy distribution to promoting renewable energy integration and accommodating varying user needs, Smart Charging embodies the convergence of sustainability, efficiency, and innovation [13].

Optimization frameworks on smart charging for EVs have seen a growing interest starting from the past decade. In [14] an heuristic method considering the State of Charge (SoC) has been proposed to minimize the charging costs. A real-time management system for the EV charging process is proposed in [15], which identifies the optimal charging periods for each EV with the aim of reducing peak load. Home charging is also worth studying due to the growing number of private EV owners, and an empirical study in [16] lead to an optimal schedule model for home charging at minimum cost. In [17], they propose an optimization model that maximizes the profit

This paper is based upon work supported by Vin University under Grant No. VUNI.2223.FT08. Furthermore, Luca Ambrosino is supported by the he FAIR - Future Artificial Intelligence Research and received funding from the European Union Next-Generation EU (PIANO NAZIONALE DI RIPRESA E RESILIENZA (PNRR) – MISSIONE 4 COMPONENTE 2, INVESTIMENTO 1.3 – D.D. 1555 11/10/2022, PE00000013). This manuscript reflects only the authors' views and opinions, neither the European Union nor the European Commission can be considered responsible for them.

ISBN:, ISSN:

of an EV aggregator in charging stations which is formulated as a mixed-integer linear programming problem.

Our research distinguishes itself by emphasizing the optimization of decision-making processes within the charging station, building a mathematical optimization model introduced in Section II, formulated as a Linear Programming problem (LP) and therefore can be solved efficiently using polynomial-time algorithms. Our objective in this paper is to assume the perspective of a charging station owner, optimizing the EVs charging process from a business point of view. We tackle the complexities of internal station management, deciding which vehicles to charge and for how long, as well as how much power to allocate to that EV, in order to help maximizing operational efficiency and reduce overall operational costs due to EVs charging process. Additionally, in Section III a Robust Optimization approach is provided for dealing with real-world uncertainties in input data such as variations in electricity prices and energy demand by EVs, encompassing the broader goal of system resilience and reliability. The main results of the cost-saving effects brought by the introduced optimization model are highlighted in Section IV, showing a huge potential of this model especially when the number of EVs grows. This innovative approach aims to make a significant contribution to the development of the field of Smart Charging both in the present and future of transportation.

## II. NOMINAL OPTIMIZATION MODEL

Some simple strategies have already been taken into account when dealing with EV's recharge, like the "First Come, First Served" (FCFS) principle [18] which ensures fairness by serving vehicles in the order they arrive. However, it doesn't consider factors like changing energy prices or demand patterns, which can lead to higher costs for charging stations. Furthermore, the risk of overloading the electricity grid must be considered. By using more flexible strategies, like adjusting prices based on demand or prioritizing urgent charging needs, stations can save money while still being fair to all users and can mitigate electricity demand peak. This helps make electric transportation more sustainable and affordable for everyone.

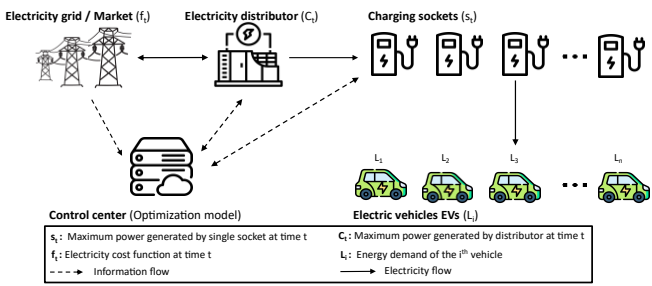

*Figure 1: Charging station structure related to our model.*

As already mentioned, this paper aims to look at the business opportunity provided by the smart charging optimization problem from the charging station's owner point of view. Figure 1 represents the scheme of the charging station as we think of it. The electricity is bought from the market and stored in the electricity distributor inside the charging station. Every charging socket in the charging station has its own power but they all share a common total power from the electricity distributor. The optimization model proposed in Eq. (1) aims to minimize the management costs of an EV charging station. These costs are associated with energy consumption, which the charging station must purchase from external suppliers following electricity market prices, assuming it does not have its own energy production, such as from solar or photovoltaic panels. Additionally, we consider that the sale of electricity from the charging station to EV customers always occurs at the same price, so that maximizing the charging station's profits is equivalent to minimizing costs.

Our optimization model deals with the power allocation to N electrical vehicles over a discrete time horizon of T steps of duration $\Delta$, considering the variable Y such that $Y_{ti}$ is the power allocated to vehicle i at time t:

$$p^* = \min_Y : \sum_{t=1}^{T} f_t \left( (1 + y_t^0) \cdot \sum_{i: t \in [a_i, d_i]} Y_{ti} \right)$$
$$s.t. \quad \sum_{t=a_i}^{d_i} Y_{ti} \geq L_i, \quad i = 1, \dots, N \quad (1)$$
$$\sum_{i=1}^{N} Y_{ti} \leq C_t, \quad t = 1, \dots, T$$
$$0 \leq Y_{ti} \leq s_t, \quad t = 1, \dots, T, i = 1, \dots, N.$$

Where:

- $f_t : \mathbb{R} \to \mathbb{R}, t = 1, \dots, T$, are convex increasing functions that encode the cost of energy at time t. In particular, if we assume that $f_t(z_t), z_t \in \mathbb{R}$ is linear $\forall t$, e.g. $f_t(z_t) = \pi_t \cdot z_t$ then the above Eq. (1) is a linear problem (LP).

- $y^0$ is a given T-dimensional-vector that corresponds to a systematic waste of energy during the charging process, for example as heat generation. We have a waste of energy if and only if the charging port is active at time t. The more energy we allocate, the more energy is waste. Theoretically, $y_t^0$ is depending on time, but for this paper we suppose it to be fixed as a small percentage as 1%.

- $L \in \mathbb{R}^N$ is the load vector, that is $L_i$ is the required energy demand $E_i$ for the i-th vehicle, divided by the time step length $\Delta$.

- $(a_i, d_i)$ are arrival and departure times for vehicle i. Of course the difference between $d_i$ and $a_i$ corresponds to the total parking time of vehicle i. Intuitively, we have $1 \leq a_i, d_i \leq T$. Address uncertainty to the EVs' arrival process rather than a scenario-based one, can be tricky and goes beyond the scope of this work.

- $C_t$ is the maximum power generated by the station at time t. $C_t$ may include a time-varying term that takes into account the self-production, coming, e.g., from solar panels. Although in general $C_t$ can be time-varying, we consider a fixed $C_t = C, t = 1, \dots, T$ when dealing with simulations, because we address the focus of the optimization on the on the time-varying electricity prices and load vectors issues.

- While $C_t$ is considered as the global charging station budget, $s_t$ instead is the maximum power generated by a single charging socket and of course is an upper bound for the variable $Y_{ti}$. For the same reason as $C_t$, we also consider a fixed $s_t = s, t = 1, \dots, T$.

The costs consist of a fixed component $y_0$, for example, due to the activation of the plug for charging or setup costs, and a variable component that naturally increases with the amount of electricity used. We assume that every charging plug has the same power, so that the cost generated by them is considered regardless to which plug is used to charge. The first constraint $\sum_{t=a_i}^{d_i} Y_{ti} \geq L_i, \quad i = 1, \ldots, N$, is a kind of demand satisfaction, explaining that the total energy delivered to each vehicle i must be greater or equal than its total energy request $L_i$. The second constraint $\sum_{i=0}^{N} Y_{ti} \leq C_t, \quad t = 1, \ldots, T$, is instead a budget constraint, assuring that at each time instant t the total power delivered by the station is never greater than the station's power capacity $C_t$. If we had wanted to emphasize the cost-saving objective more in the optimization model, we would have avoided imposing a direct budget constraint, opting instead to incorporate a penalty into the objective function for exceeding the budget. In our context, it's important to highlight that the goal of the model is also to promote energy sustainability by preventing grid overload, despite the potential for greater cost savings from overloading. This balance is crucial to ensure the long-term stability and efficiency of the energy system, thereby contributing to the overall sustainability of the energy infrastructure. Finally, $0 \leq Y_{ti} \leq s_t, \quad t = 1, \ldots, T, i = 1, \ldots, N$, states that each variable $Y_{ti}$ must obviously be non-negative, since it represents power directionally allocated to vehicle i at time t, and has $s_t$ as upper bound.

In order to simplify the representation and the simulations code, we use a binary matrix A of size $T \times N$ to encode arrival and departure times $(a_i, d_i)$. Each element $A_{ti}$ equals 1 if t falls into the interval $[a_i, d_i]$ for vehicle i, 0 otherwise. This matrix allows us to condense the information regarding the arrival and departure times of each vehicle into a more compact format. With this matrix representation, we can reformulate the optimization problem stated above in a more concise matrix form useful for Matlab simulations:

$$\min_Y: F\big(diag(AY^T) + y^0 \cdot diag(AY^T)\big)$$
$$s.t.: \ diag(A^TY) \geq L,$$
$$diag(AY^T) \leq C, \qquad (2)$$
$$Y \geq 0,$$
$$Y \leq s.$$

Where:

- $F: \mathbb{R}^T \to \mathbb{R}$ aggregates all the cost functions $f_t$, that is $F(z) := \sum_{t=1}^{T} f_t(z_t), \quad z \in \mathbb{R}^T$.
- The diagonal elements of the N×N-dimensional matrix $A^TY$, say $(A^TY)_{ii}$, are the total power absorbed by the vehicle i during its whole charging.
- The diagonal elements of the T×T-dimensional matrix $AY^T$, say $(AY^T)_{tt}$, are the total power available at time t.
- L and C are respectively the vectors related to the power requirement and the capacity (i.e. energy available). Finally, $s$ is the vector related to the single socket maximum power generated.

III. ROBUST OPTIMIZATION OVERVIEW

The sources of uncertainty considered in this optimization model are related to electricity prices $\pi_t$ and to load vectors $L_i$, i.e. energy required by EVs. The electricity price fluctuates not only periodically with time (during a day of operation there are peak-price periods and low-price periods, e.g., at night) but also in response to unforeseeable network and real-time electricity market conditions, which can be difficult to forecast. Concerning the energy demand from EVs, we do not know exactly in advance how much energy each vehicle would need to be charged. Several factors affect this uncertainty and a huge problem may arise during peak hours if there are many EVs in the charging station needing a great amount of charge. These uncertainties can be characterized by the fact that the objective function and constraints of the problem are dependent on additional parameters $u \in \mathcal{U} \subset \mathbb{R}^d$, which represent the uncertainties. The functions $f_i(x, u)$, which can be either objective function or constraint, represent the uncertain components of the problem. Depending on the approach chosen to define and manage uncertainty, various methods for optimization under uncertainty are available.

In this paper, the two robust optimization methods used to deal with uncertainties in electricity prices and in load vector are respectively the norm-bounded uncertainty in dimension 2 and the interval uncertainty [19]. In formula, a norm-bounded uncertainty set (a.k.a. norm ball) is a set of the form:

$$\mathcal{U}: \{u: \|u - \hat{u}\|_p \leq r\} = \{r = \hat{u} + rz: \|z\|_p\} \leq 1 \qquad (3)$$

Where $\hat{u}$ is the center of the ball, $r$ is its radius, and $p$ denotes the type of norm, typically $p = 1, 2, \infty$. In particular, in the common case when $p = 2$ the norm ball is a hypershpere and the problem is a second order cone programming (SOCP) problem. For our specific smart charging problem in Eq. (2), we assume that $\pi$, i.e. the time-varying vector encoding electricity prices, is only known up to a sphere: $\|\pi - \hat{\pi}\|_2 \leq r$, where the "nominal" cost $\hat{\pi} \in \mathbb{R}^T$ and the uncertainty level $r \geq 0$ are both known. The robust counterpart to the model in Eq. (2) with linear prices $\pi$ is the SOCP:

$$\min_Y: \pi \cdot \big(y^0 + diag(AY^T)\big) + r \cdot \|y^0 + diag(AY^T)\|_2$$
$$s.t.: \ diag(A^TY) \geq L,$$
$$diag(AY^T) \leq C, \qquad (4)$$
$$Y \geq 0,$$
$$Y \leq s.$$

We observe that when $r = 0$, we recover the nominal problem with fixed costs, because the uncertainty ball has radius 0, while for $r$ large, the solution drives toward the problem with quadratic costs instead of linear.

In practice, due to uncertainty on the initial charge of vehicles, also the load vector L is uncertain. Assume for example that the load vector L is only known to satisfy $L \in [\underline{L}, \overline{L}]$, where the lower and upper bounds $\underline{L}$ and $\overline{L}$ are known somehow. The robust counterpart has a trivial expression, in the same form as Eq. (4), but with values of L replaced with their worst-case (largest) values $\overline{L}$, which means:

$$\min_Y: \ \pi \cdot \big(diag(AY^T) + y^0 \cdot diag(AY^T)\big) +$$
$$+ r \cdot \|diag(AY^T) + y^0 \cdot diag(AY^T)\|_2$$
$$s.t.: \ diag(A^TY) \geq \overline{L}, \qquad (5)$$
$$diag(AY^T) \leq C,$$
$$Y \geq 0,$$
$$Y \leq s.$$

In some practice applications, if the infrastructure of the charging station allows to measure, at charging time, the actual desired load from EVs, then an affine recourse strategy may improve the performance of robust optimization.

IV. SIMULATIONS AND RESULTS

After the optimization model theoretical introduction, we want to show the benefits that this smart charging approach can give to one charging station's finances. With this aim, we perform numerical simulations comparing our model in Eq. (1) to a more trivial charging decision algorithm like the First Come First Served (FCFS) process. This standard algorithm states that, whenever an EV enters the charging station, it is put into charge if the total amount of energy deliverable by the charging station allows it. In formula, the power allocated to vehicle $i$ at time $t$, say $X_{ti}$, is:

$$X_{ti} = min([s_t, L_i^{res}, C_t^{res}]) \qquad (7)$$

If the charging station is full when an EV arrives, it must wait until another EV completes its recharge. In fact, when the total power deliverable by the charging station $C_t$ is already in use for other Evs, we will have $C_t^{res} = 0$ in Eq. (7) forcing into $X_{ti} = 0$. $L_i^{res}$ is instead the residual energy required by vehicle $i$ considering the charging already done in the previous time steps, so that if an EV needs a small amount of energy, the algorithm gives it only the amount needed in order to avoid waste of energy and useless costs. The decision process in Eq. (7) is a standard procedure currently adopted by many charging stations which do not use smart charging methods. It is a basic method which already takes care about the energy waste, therefore is better than other simpler algorithms. The weakness of this method is that it does not consider the changing electricity prices during the different time instants $t$, instead it aims to charge the EVs as soon as possible. This can be a reasonable decision, but if we know that an EV stays parked in the station for a long time, a smart charging method that considers time-varying electricity prices can improve the performances and optimize the cost-saving.

The operational cost minimization is obviously related also to a lower use of electricity and leads to a greener economy. To show the benefit of our approach, we run MATLAB simulations comparing the optimization model in Eq. (2) with the more trivial algorithm just introduced in Eq. (7). The input data considered are taken from the Caltech university charging station in California (USA) [20], which provides a common and open-source dataset useful to conduct analysis on EVs. The specific dataset that we used to generate in our simulations both the load vector $L$ and the EVs arrival scenarios $A$ contains records for 861 days, starting from *25-Apr-2018*. The final day recorded in our Caltech.csv dataset is *12-Apr-2021*. The dates within the dataset are not always consecutive, as there are 223 missing dates. Most of the missing dates are after *19-Mar-2020*, which marks the beginning of the lockdown in California due to the COVID-19 pandemic [21]. Other missing days may be attributed to days when no electric vehicles visited the charging station, especially on Sundays. Despite the gaps, we can consider each recorded day as a separate scenario, each providing its own set of input data for the model. This means that each day, regardless of continuity, offers a unique dataset that can be independently analyzed and modeled.

The dataset contains records of EVs' arrival and departure time, charging duration and total energy absorbed in kWh. From this data, we can easily compute the input matrix $A$ and the load vector $L$ of energy required by EVs for each day in kWh, considering a fixed time horizon $T = 24$ hours, i.e. 24 optimization time steps. The other parameters of the model have been fixed: $y_t^0$ even if theoretically can vary in time, has been set to a reasonable small value; $C_t$ and $s_t$ can also depend on time, but we chose to refer to the exact values of Caltech charging station [22]. It is not easy to extract electricity prices data from California, since we need hour-by-hour prices information. The input vector $\pi$ used comes from the open source European electricity prices [23], specifically from Italy, considering the same day as matrix $A$ and vector $L$ refer to, in order to have coherent scenarios. Be aware that, even if the electricity prices do not correspond to the same geographical area of the EVs data, our model is consistent and general enough to be applied to different electricity prices. In fact, different electricity prices schemes only affect the percentage of cost-saving, and do not undermine the model efficiency. All input data are summarized in Table 1: Input values for simulations.

*Table 1: Input values for simulations*

| Input data | Value |
|---|---|
| $A, L$ | Estimated from Caltech data |
| $y_t^0$ | $0.01, \forall t$ |
| $C_t$ | $300\ kW, \forall t$ |
| $s_t$ | $7\ kW, \forall t$ |
| $\pi$ | Italian electricity prices in € |

Already analyzing the simulation performed on the first day of the dataset, i.e. *25-Apr-2018*, it is clear that the time-aware approach with respect to electricity prices benefits a lot, compared to the trivial algorithm decision process in Eq. (7). In fact, the huge cost-saving on this specific day is 28.98%.

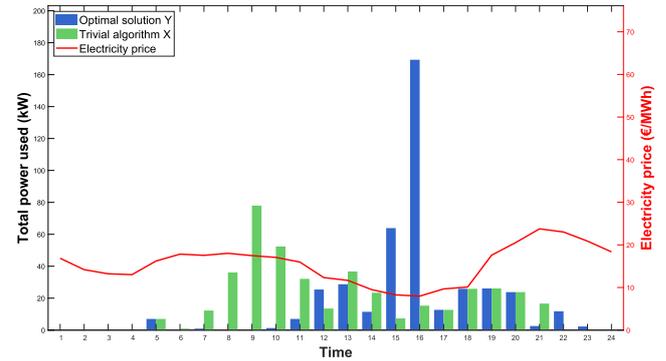

*Figure 2: Total power allocated to EVs during one single day (25-Apr-2018), divided by hours.*

Figure 2: Total power allocated to EVs during one single day (25-Apr-2018), divided by hours. shows the total power allocated by the two methods over the whole charging station, for each time step. We can notice that while the trivial algorithm allocate power to charge EVs as soon as they enter the charging station, our optimization model aims to charge as much as possible when the electricity price is lower during the day.

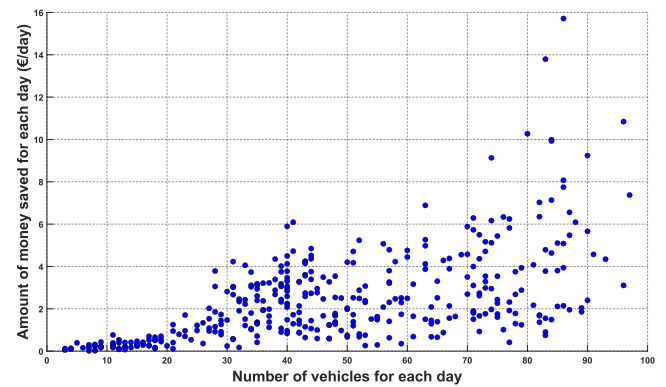

*Figure 3: Total money saved for each day, over 365 scenarios simulations, starting from 25-Apr-2018.*

It is easy to understand that applying this reasoning to a bigger number of EVs, manages to achieve an higher cost-saving. The scatter plot in Figure 3 shows the positive correlation between the number of EVs entering the charging station during the whole day and the amount of money saved on that day, due to our model's smart allocation of the recharging power. Each blue dot in Figure 3 represents one of the 365 days over the simulation horizon. The cost-saving percentage obtained by applying our optimization model every day for 365 days scenarios is 8.78% on average, with a maximum cost-saving of 55.25% and a minimum of 0.66%, obtained on 31-Dec-2018 when only 7 EVs visited the charging station due to New Year's Eve.

The scatter plot in Figure 3 implicitly suggests a huge potential of smart charging for the future, when the number of Electric Vehicles will increase in order to obtain a greener economy for a long-term sustainable world. This smart charging optimization model fits well the long-term needs because we can save always more money when we apply this power allocation decision process over a longer time period, as shown in the area plot in Figure 4. A 2-year long simulation period, i.e. 730 different scenarios, shows the potential of our smart charging mathematical model when it is continuously applied to reduce the operational costs of a charging station. The optimized cost are always lower than the cost obtained by applying the trivial decision algorithm for EVs charging, and this gap is growing in time also due to the higher number of EVs on the market. The average cost-saving percentage over 730 scenarios is 8.95%.

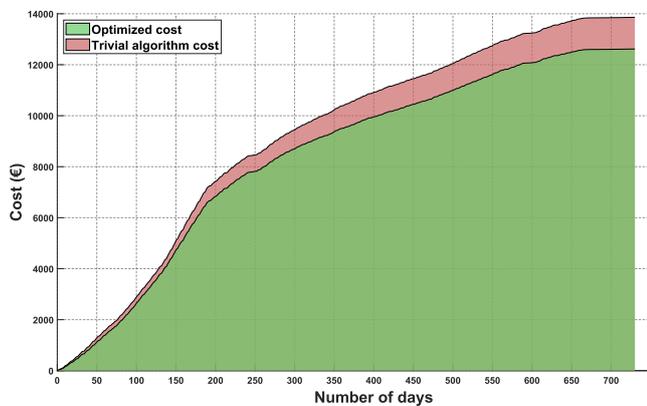

***Figure 4:*** *Total operational cost in time, over a 730 scenarios period starting from 25-Apr-2018.*

In Figure 4 we can notice that the cost is increasing quite fast in the first half of the simulation time period, while it becomes more flat in the second half. This is surely an effect of fewer EVs coming to the charging station after the outbreak of COVID-19 pandemic. In fact, not only after the *19-Mar-2020* there are a lot of missing days in the dataset, but there are reasonably also many days with very few EVs. More in details, the average number of EVs among the 861 scenarios is 30.1. Another computational simulation has been performed considering only the 455 scenarios with number of EVs $N \geq 30$, obtaining an average cost-saving percentage of 9.88%. A summary of the simulation results is exposed in Table 2, comparing the different costs obtained applying either the trivial decision algorithm or the smart charging optimization model, with same input parameters as in Table 1 but different settings related to number of scenarios $T$ and to filtered condition on the number $N$ of EVs in each scenario. Here the *Trivial Cost* and the *Optimized Cost* are computed by summing up all the daily operational costs over the number of scenarios $T$. The numbers in Table 2 indicate once again that the smart charging optimization model has a huge potential in cost-saving during the days with a lot of EVs in the charging station. In fact, in the dataset there are 574 scenarios out of 861 where the number of electric vehicles (EVs) is at least 10, accounting for exactly two-thirds of the dataset. However, these 574 scenarios generate almost 97% of the total operational costs either we consider the trivial algorithm or the mathematical optimization model.

***Table 2:*** *Simulations results summary.*

| T | EVs Number | Trivial Cost | Optimized Cost | Average Cost-Saving |
|---|---|---|---|---|
| 100 | N>0 | 2890.03€ | 2642.52€ | 8.56% |
| 365 | N>0 | 10454.64€ | 9554.42€ | 8.61% |
| 861 | N>0 | 14033.47€ | 12780.11€ | 8.93% |
| 574 | N≥10 | 13599.24€ | 12370.08€ | 9.04% |
| 455 | N≥30 | 12579.39€ | 11430.53€ | 9.13% |

## V. CONCLUSIONS AND FUTURE WORKS

Smart Charging is a research field that has recently attracted the attention of the world and is developing day by day. The EVs market is still at the beginning of its expansion, but it is only a matter of time that it is going to exponentially grow, since the climate crisis is striking more and more, and many governments have already adopted some preventive measures that will lead to an always higher number of EVs in the streets. Of course, charging stations already adopt decision algorithms to determine a strategy to follow when many EVs approach the station. Anyway, not always the decision made turns out to be optimal. In fact, we have shown that following a Smart Charging optimization model significantly reduces the charging station operational cost, while keeping the same service level and customer satisfaction.

The mathematical optimization model presented in this paper utilizes a time-aware strategy to efficiently allocate the energy needed for each EV to recharge, intelligently selecting the time interval and vehicle to charge in order to minimize the operational costs associated with purchasing electricity from the power grid by the charging station. This smart charging model indeed ensures an average 9% cost saving compared to a standard algorithm when tested across 730 days, i.e. 730 different scenarios related to vehicle arrivals, EVs energy requirement, and electricity costs. Moreover, this model, allows for good flexibility to be used within the realm of robust optimization, aiming to find an optimal smart charging policy even when the input parameters are uncertain.

Although our model reach interesting results in optimization, there are some aspects that deserve to be taken into account as potential improvements. The objective function of the proposed optimization model focuses on the cost-saving without considering the profitability of the charging station. If the electricity is sold to the customer at a fixed price, then this approach can be considered optimal. However, challenges arise when the prices of electricity vary over time or if discounts are offered for extended charging periods. In that case the potential of our smart charging optimization model would significantly grow, but it could be convenient to re-formulate the optimization problem under a profit maximization point of view for the objective function. This idea can easily fit a company needs in the real-world, but these non-linearities are going to significantly complicate the model and its computational resolution.

Another consideration about the proposed optimization model is the assumption that all charging sockets have uniform power output. However, it's worth noting that many charging stations currently offer various types of sockets with different charging capabilities, such as slow, standard, or fast charging. These distinctions not only affect the charging speed but can also influence pricing, introducing additional complexity to the optimization process. Our model already incorporates factors like electricity cost fluctuations over time and manages to efficiently schedule EV recharges based on their parking durations, as outlined in matrix A. Should the charging station feature different socket types, exploring the interplay between EV energy requirements and parking durations could potentially yield even more optimal charging schedules. As already mentioned, the EV market is likely to grow in the next few years, and when this happens it could be useful to consider adding an upper bound constraint on the cardinality of the optimization problem solution. Some constraint relaxation techniques may help keeping the computational problem at an acceptable difficulty despite a non-linear cardinality constraint.

Considering what we have observed and the few assumptions made to deeply analyze the performances of this mathematical model, it is important to emphasize that this smart charging model aspires to be a benchmark for other optimization strategies related to Smart Charging for EVs. Certainly, developing more complex models based on this one's idea will be a central goal for this research field's future, which aims to reconcile humanity's daily practical transportation needs with the global shared objective of shifting the world economy towards long-term energy sustainability to safeguard the health of our planet.